\documentclass[11pt]{amsart}
\usepackage[latin1]{inputenc}

\usepackage{amsmath}
\usepackage{amsfonts}
\usepackage{color}
\usepackage{amssymb}
\usepackage{amsthm}
\usepackage{graphicx,epsfig}
\usepackage{amscd}
\newcommand{\nil}{\mathrm{Nil}_3}

\newcommand{\RA}{\R^2\rtimes _A \R}

\newcommand{\R}{\mathbb{R}}

\newcommand{\N}{\mathbb{N}}

\newcommand{\rth}{\mathbb{R}^3}

\newcommand{\cS}{{\mathcal S}}

\newcommand{\cA}{{\mathcal A}}

\newcommand{\cB}{{\mathcal B}}
\newcommand{\cC}{{\mathcal C}}

\newcommand{\cP}{{\mathcal P}}

\newcommand{\ben}{\begin{enumerate}}
\newcommand{\bit}{\begin{itemize}}
\newcommand{\een}{\end{enumerate}}
\newcommand{\eit}{\end{itemize}}

\newcommand{\ve}{\varepsilon}

\newcommand{\ed}{\end{document}}

\newtheorem{thm}{Theorem}[section]

\newtheorem{theorem}{Theorem}[section]

\newtheorem{lemma}[thm]{Lemma}

\newtheorem{definition}[theorem]{Definition}

\newtheorem{remark}[theorem]{Remark}

\newtheorem{exa}[thm]{Example}
\newtheorem{corollary}[theorem]{Corollary}
\newtheorem{conjecture}[theorem]{Conjecture}

\newtheorem{cla}[thm]{Claim}

\newcommand{\wt}{\widetilde}
\newcommand{\HH}{\mathbb{H}}
\def\esf{\mathbb{S}}

\begin{document}

\title[Half-space theorems and the embedded Calabi-Yau problem]{Half-space theorems and the embedded Calabi-Yau problem in Lie groups}
\author{Beno\^\i t Daniel}\thanks{The first author is partially supported by Minist\`ere des Affaires
\'etrang\`eres et europ\'eennes (France), Partenariat Hubert Curien STAR,
projet n${}^{\circ}$ 23889UG}
\author{William H. Meeks III}
\thanks{This material is based upon work
 for the NSF under Award No. DMS-1004003 by the second author. Any opinions,
 findings, and conclusions or recommendations expressed in
 this publication are those of the authors and do not necessarily
 reflect the views of the NSF}
\author{Harold Rosenberg}

\date{\today}


\address{Universit\'e Paris-Est Cr\'eteil, D\'epartement de Math\'ematiques, UFR des Sciences et Technologies, 61 avenue du G\'en\'eral de Gaulle, 94010 Cr\'eteil cedex, FRANCE}
\email{daniel@univ-paris12.fr}
\address{Mathematics Department, University of Massachusetts, Amherst, MA 01003, USA}
\email{profmeeks@gmail.com}
\address{Instituto Nacional de Matem\'atica Pura e Aplicada (IMPA), Estrada Dona Castorina 110, 22460-320 Rio de Janeiro - RJ, BRAZIL}
\email{rosen@impa.br}

\begin{abstract}
We study the embedded Calabi-Yau problem for complete embedded constant mean curvature surfaces of finite topology or of positive injectivity radius in  a simply-connected three-dimensional Lie group $X$ endowed with a left-invariant Riemannian metric.
We first prove a half-space theorem for constant mean curvature surfaces. This half-space theorem applies to certain properly immersed constant mean curvature surfaces of $X$ contained in the complements of normal $\R^2$ subgroups $F$ of $X$.
In the case $X$ is a unimodular Lie group, our results imply that every  minimal surface in $X-F$ that is properly immersed in $X$ is  a left translate of $F$
and that every complete embedded minimal surface of finite topology or of positive injectivity radius in  $X-F$ is also a left translate of $F$.
\end{abstract}


\maketitle

\section{Introduction.} \label{sec:intro}

A natural question in the global theory of minimal surfaces,
first raised by Calabi in 1965~\cite{ca4} and later revisited by
Yau~\cite{yau1, yau2}, asks whether or not there exists a complete
immersed minimal surface $\Sigma$ in a bounded domain of $\rth$, or more generally, the question asks: {\em If $\Sigma$  is contained in a half-space of $\rth$, then is it a plane parallel to the boundary of the half-space?} For complete {\em immersed} minimal surfaces these questions were answered by the existence results of Jorge and Xavier~\cite{jx1} and by Nadirashvili~\cite{na1}. Closely related to the Calabi-Yau problem in Euclidean $3$-space is the half-space theorem by Hoffman and Meeks~\cite{HM}: {\em If $\cS$ is a properly immersed minimal surface in $\R^3$ that lies on one side of some plane $\cP$, then $\cS$ is a plane parallel to $\cP$.} A fundamental result  on the embedded Calabi-Yau problem for complete embedded minimal surfaces in $\rth$ was recently given by Colding and Minicozzi~\cite{cm35}: {\em A complete embedded minimal surface $\cS$ of finite topology in $\rth$ is proper. In particular, if such an $\cS$ lies in a half-space of $\rth$, then it is a plane.}  Meeks and Rosenberg~\cite{mr13} then generalized this last theorem by proving that any complete embedded minimal surface of positive injectivity radius in $\rth$ is proper; since a complete immersed minimal surface of finite topology in $\R^3$ always has positive injectivity radius, this last result generalizes the previously stated  theorem of Colding and Minicozzi.

In this paper we generalize our work in~\cite{dmr1} and the above classical results about minimal surfaces to simply-connected homogeneous manifolds. We recall that every simply-connected homogeneous three-manifold $X$ is either isometric
to $\esf^2\times \R$ with a constant curvature metric on $\esf^2$ or it is a metric Lie group, where by metric Lie group we mean  a Lie group $X$ equipped with a left-invariant Riemannian metric. We generalize our work in~\cite{dmr1} by  showing that whenever $F$ is a normal $\R^2$ subgroup of a three-dimensional simply-connected unimodular Lie group $X$   and $\cS$ is a properly immersed minimal surface in $X$ that lies on one side of $F$, then $\cS$ is a left translate of $F$; see item~\ref{it1} of Theorem~\ref{main} below and the sentence just before the statement of that theorem. We then go on to generalize the above stated results of Colding-Minicozzi and of Meeks-Rosenberg to prove that if $\cS\subset X$ is a complete embedded minimal surface of finite topology or of positive injectivity radius which lies on one side of $F$, then $\cS$ cannot have $F$, or a left translate of $F$, in its closure, and so, by the proof of our half-space theorem, $\cS$ is a left translate of $F$.

We  call a connected surface $\cS \subset X$ of constant mean curvature $H$ an {\em $H$-surface}; after appropriately orienting $\cS$, we will always assume that  $H\geq 0$. Suppose that $F$ is a  normal $\R^2$ subgroup in a simply-connected three-dimensional metric Lie group $X$. Theorem~\ref{main} below gives conditions for certain properly immersed $H$-surfaces of $X$ contained in  $X-F$  to be left translates of $F$. We note that whenever a metric Lie group $X$ admits such a normal $\R^2$ subgroup, then $X$ can be expressed as a semidirect product $\RA$ of this normal subgroup with $\R$, where the related homomorphism $\sigma \colon\R\, \to \mbox{Aut}(\R^2)=\mathrm{GL}(2,\R)$ is given by $\sigma(t)=e^{tA}$ for some $2\times2$ matrix $A$. For this last result see~\cite{mpe11,mil2} and see~\cite{mpe11} for how one obtains a left-invariant metric from the matrix $A$. For the sake of being self-contained, in Section~\ref{pre} we  review the main constructions  related to $\RA$ and other useful formulas from~\cite{mpe11} that are needed in the proof of the next theorem. After possibly replacing $A$ by $-A$, which changes the orientation of $X=\RA$, {\bf we will always assume} trace$(A)\geq 0$.

Recall that every simply-connected homogenous three-manifold $X$ is either isometric to $\esf^2\times \R$ with a constant curvature metric on $\esf^2$ or it is a metric Lie group.  We define:
\begin{center} $H_0(X)$  is the infimum of the mean curvatures of immersed  compact constant mean curvature surfaces in $X$ with $H\geq 0$.\end{center}
By results in~\cite{mmpr1}, when $X$ has the form $\RA$ , then $H_0(X)=\frac12$trace$(A)$ and for $X=\mathrm{SU}(2)$ with a left-invariant metric or for  $X=\esf^2\times \R$, then $H_0(X)=0$.  The value $H_0(X)$ is only not known when the group structure of $X$ is that of $\wt{\mathrm{PSL}}(2,\R)$ and the metric is a general left-invariant metric.

For the purpose of understanding and applying Theorem~\ref{main} stated below, we remark that a normal $\R^2$ subgroup $F$ of a metric Lie group $X$ has an intrinsic flat metric, has constant mean curvature $H_0(X)$; furthermore,  $F$ is minimal if and only if $X$ is a unimodular Lie group; see~\cite{mmpr1,mpe11} or the appendix in~\cite{mpe7} for the proof of this last fact about unimodular groups.

\begin{theorem} \label{main}
Suppose $F$ is a left translate of a normal  $\R^2$ subgroup of nonnegative constant mean curvature $H_0(X)$ in a metric Lie group $X$.  Then:
\ben
\item \label{it1} If $H_0(X)=0$, then every properly immersed minimal surface in $X$ that is disjoint from $F$ is a left translate of $F$.
\item  \label{it2} If $H_0(X)>0$ and $H\in [0, H_0(X)]$, then every properly immersed $H$-surface $\cS$ in $X$ that lies on the mean convex side  of $F$ is a left translate of $F$. In particular,  the mean curvature of such an $\cS$ is $H_0(X)$.
\item  \label{it3} If  $H_0(X)>0$ and $H\in [H_0(X), \infty)$, then  every properly {\em embedded} $H$-surface $\cS$ in $X$ such that $F$ lies on the mean convex side  of $\cS$ is a left translate of $F$. In particular,  the mean curvature of such an $\cS$ is $H_0(X)$.
\een
\end{theorem}

Previously, half-space theorems in  three-dimensional metric Lie groups $X$ have played an important theoretical role in the proofs of deeper geometric and analytic results; for example, see the papers \cite{dh,dmr1,hrs,espinoza,roro} for examples and applications to the Bernstein problem and related results.  In relationship to the half-space theorem in~\cite{espinoza}, we recall that there exist left-invariant metrics on the Lie group $\wt{\mathrm{PSL}}(2,\R)$ which have four dimensional isometry groups and which, as homogeneous three-manifolds, are isometric (but non-isomorphic as groups) to the metric Lie groups $\RA$, where the matrix $A$ has the form:
$A=\left(
\begin{array}{cc}
a & 0\\
c & 0\\
 \end{array}\right), \;\mbox{for}\; a>0 \; \mbox{for}\; c\neq 0$.  Note that as $a,c$ vary the Lie group structures for the metric Lie groups $\RA$  are the same, only the metrics are changing.

It is well-known that for each such $A$, there exists a Riemannian submersion $\Pi \colon \RA \to \HH^2(\kappa(A))$ with some constant bundle curvature $\tau(A)$ and where $\HH^2(\kappa(A))$ denotes the hyperbolic plane equipped with a metric of constant negative curvature $\kappa(A)$; see item~2 of Theorem~2.13 in~\cite{mpe11} for this result.
Under the isometric correspondence between this metric Lie group $\RA$ and the related metric Lie group $\wt{\mathrm{PSL}}(2,\R)$ with four dimensional isometry group, the normal subgroup $\R^2 \rtimes_A \{0\}$ corresponds to a horocylinder which is the inverse image under $\Pi$ of a particular horocycle in the base space $\HH^2(\kappa(A))$.  Thus, we can apply Theorem~\ref{main}  even in the case of the unimodular group $\wt{\mathrm{PSL}}(2,\R)$ which has no normal $\R^2$ subgroups and  the above theorem still give a useful half-space theorem; this half-space theorem  was recently proved by Espinoza~Pe\~nafiel~\cite{espinoza} using different methods. When $c$ is allowed to be $0$ in the above matrix $A$, then $\RA$ is isometric to $\HH^2(-a^2)\times \R$, where $\HH^2(-a^2)$ denotes the hyperbolic plane with constant Gaussian curvature $-a^2$ and we recover the halfspace theorem in~\cite{hrs} for $H$-surfaces which lie on one side of a horocylinder in $\HH^2(-a^2)\times \R$.

Our second theorem concerns the embedded Calabi-Yau problem.
The next related conjecture is known to hold for $X=\esf^2\times \R$, $X=\rth$ and for any metric Lie group $X$ with a constant curvature metric (see~\cite{cm35,mpr19,mr13,mt11}). It is also natural to conjecture properness in the next conjecture with ``finite topology" replaced by ``positive injectivity radius".

\begin{conjecture}[Embedded Calabi-Yau problem for homogeneous 3-manifolds] \label{con:CY}
A complete embedded $H$-surface $\cS$ of finite topology in a simply-connected homogeneous $3$-manifold $X$ is proper if $H\geq H_0(X)$. In particular, if $X$ is {\em SU}$(2)$ with a left-invariant metric, then $\cS$ is compact.
\end{conjecture}

Recently Meeks and Tinaglia~\cite{mt11} proved Conjecture~\ref{con:CY} holds in the special case where $X$ is the hyperbolic three-space $\mathbb{H}^3 =\R^2 \rtimes_A \R$, where $A$ is the $2\times 2$ identity matrix. This proof uses a halfspace theorem in $\mathbb{H}^3$ by
Rodriguez and Rosenberg~\cite{roro}, which can be viewed as a special case of the halfspace theorem described in Theorem~\ref{main}. In order to prove Conjecture~\ref{con:CY} holds in this special case, Meeks and Tinaglia applied this halfspace theorem of Rodriguez and Rosenberg to demonstrate: A necessary and sufficient condition for a complete, connected  embedded $H$-surface in $\mathbb{H}^3$ with $H\geq 1$ to be proper is that it have injectivity radius function bounded away from zero on compact subsets of  $\mathbb{H}^3$, a property that holds for complete embedded finite topology $H$-surfaces in any metric Lie group (see Section~\ref{sec:CY} and Theorem~\ref{mp} for further discussion).

The following theorem is related to the above conjecture for the value $H=H_0(X)$.  This relationship will become clear in the proof of Theorem~\ref{CY}.

\begin{theorem} \label{CY}
Suppose $F$ is a  left translate of a normal  $\R^2$ subgroup in a three-dimensional, simply-connected metric Lie group $X$ and  $\cS \subset X$ is a complete embedded $H$-surface with $H \in[0,H_0(X)]$ such that $\cS$ is contained on the mean convex side of $F$ (both sides are mean convex if $F$ is minimal).  \ben \item If $\cS$ has finite topology or if it has positive injectivity radius, then $\cS$  is a left translate of $F$.
  \item If $H \in(0,H_0(X)]$ and $\cS$ has  injectivity radius function bounded away from zero on compact domains of $X$, then $\cS$  is a left translate of $F$.
      \een
\end{theorem}

An example of the relationship of the embedded Calabi-Yau problem to  Theorem~\ref{CY} and the next theorem can be seen as follows. Recall that the Heisenberg group with its standard metric is $\nil=\RA$, where $A=\left(
\begin{array}{cc}
0 & 1\\
0 & 0\\
 \end{array}\right)$.   If $\cS$ is a complete embedded $H$-surface in $\nil$ of positive injectivity radius or of finite topology which is not proper, then, by the results in~\cite{mpr18,mt1}, its closure $\overline{\cS}$ contains a complete embedded stable $H$-surface $\Sigma \subset (\nil - \cS)$. It is conjectured that such a stable minimal surface $\Sigma$  would be a entire minimal graph with respect to a certain natural Riemannian submersion $\Pi \colon \nil \to \R^2$ or, after an isometry of $\nil$, would be $\R^2 \rtimes_A \{0\}$. Hence, by Theorem~\ref{CY} and Corollary~\ref{cor:CY} below, after an isometry of $\nil$, $\cS$ would be a translation of $\Sigma$, and thus, $\cS$ is in fact proper in $\nil$.


For the statement of the next theorem, recall that a surface $F$ in $\nil$ is called an {\em entire graph} if $\Pi|_{F}\colon F\to \R^2$ is a diffeomorphism. Note that in this case the fibers of this fibration are considered to be vertical (vertical is different from the third coordinate direction when we view $\nil$ as being $\RA$).

\begin{theorem} \label{CY2}
Suppose $F$ is an entire minimal graph in the metric Lie group $\nil$ and  $\cS \subset X-F$ is a complete embedded minimal surface.  If $\cS$ has positive injectivity radius, then $\cS$ is proper in $\nil$.
\end{theorem}

\begin{corollary} \label{cor:CY}
Under the hypotheses of Theorem~\ref{CY2}, $\cS$ is a vertical translate of $F$.
\end{corollary}

This corollary follows from  Theorem~1.4 in~\cite{dmr1}.
Since a complete embedded minimal surface in $\nil$ of finite topology has positive injectivity radius \cite{mpe7},  we can replace "positive injectivity radius" by "finite topology" in this theorem (see Section \ref{sec:CY} for details).

The paper is organized as follows. In Section~\ref{pre} we recall the basic description and related formulae for metric Lie groups of the form $\RA$; this material is taken from the paper~\cite{mpe11}, where the reader can find further discussion and details. In Section~\ref{sec:halfspace},  Theorem~\ref{main} is proved and in Section~\ref{sec:CY}, we prove Theorems~\ref{CY} and~\ref{CY2}.

We remark that Mazet~\cite{maz3} has independently proven a result, which together with our discussion of Theorem~\ref{main}, implies the conclusions of that theorem hold; in fact, his main theorem is a general half-space type theorem for properly immersed $H$-surfaces in certain Riemannian three-manifolds which holds in greater generality than the case where the manifold is a metric Lie group.

\vspace{.2cm}

\noindent {\bf Acknowledgments.} The first two authors would like to thank the Korea Institute for Advanced Study (KIAS) where part of this work was completed.

\section{Preliminaries.}\label{pre}
In this section we recall some useful formulas for three-dimensional Lie groups of the form $X=\RA$ that we   need in our proofs. This discussion and the related formulas are taken from~\cite{mpe11} and we refer the reader to that paper and to~\cite{mil2} for further details.

Generalizing direct products, a semidirect product is a particular
way of cooking up a group from two subgroups, one of which is a
normal subgroup. In our case, the normal subgroup will be
isomorphic to $\R^2$  and the
other factor  will be isomorphic to $\R $. As a set, this semidirect
product is nothing but the cartesian product of $\R^2$ and $\R$, but the group
operation is different. The way of gluing different copies of $\R^2$ is
by means of a 1-parameter subgroup $\varphi \colon \R \to
\mbox{Aut}(\R^2)$ of the automorphism group of $\R^2$:
\[
\begin{array}{rccc}
\varphi (z)=\varphi _z\colon &\R^2&\to &\R^2\\
& {\bf p }&\mapsto &\varphi _z({\bf p}).
 \end{array}
\]
Hence the group operation of the semidirect product $\R^2\rtimes
_{\varphi }\R$ is given by
\begin{equation}
\label{eq:5}
 ({\bf p}_1,z_1)*({\bf p}_2,z_2)=({\bf p}_1 + \varphi
_{z_1} ({\bf p}_2),z_1+z_2),
\end{equation}
where  $\,+\,$  denotes the group operation in both $\R^2$ and $\R$.
Since $\mbox{Aut}(\R^2)=\mbox{Gl}(2,\R)$, $\varphi $ is given by exponentiating a given
matrix $A\in {\mathcal M}_2(\R )$, i.e., $\varphi _z({\bf p})=
 e^{zA}{\bf p}$, and we   denote the corresponding group by
$\RA$.

Our first goal is to describe the left and right-invariant vector fields on
a semidirect product $\R^2\rtimes _A\R $ for any matrix
 $A=\left(
\begin{array}{cr}
a & b \\
c & d \end{array}\right) $. To carry out  this goal, first choose coordinates $(x_1,{x_2})\in\R^2$, ${x_3}\in \R $. Then
$\partial _{x_1}$, $\partial _{x_2},\partial
_{x_3}$ is a parallelization of $X=\R^2\rtimes _A\R $. Taking
derivatives at $t=0$ in the expression (\ref{eq:5}) of the left
multiplication by $({\bf p_1},{x_3}_1)=(t,0,0)\in X$ (resp. by $(0,t,0),
(0,0,t)$), we obtain the following basis $\{ F_1,F_2,F_3\} $ of the
right-invariant vector fields on $X$:
\begin{equation}
\label{eq:66}
 F_1=\partial _{x_1},\quad F_2=\partial _{x_2},\quad F_3({x_1},{x_2},{x_3})=
(a{x_1}+b{x_2})\partial _{x_1}+(c{x_1}+d{x_2})\partial _{x_2}+\partial _{x_3}.
\end{equation}
Analogously, if we take derivatives at $t=0$ in the right
multiplication by $({\bf p_2},{x_3}_2)=(t,0,0)\in X$ (resp. by $(0,1,0),
(0,0,1)$), we obtain the following basis $\{ E_1,E_2,E_3\} $ of the
Lie algebra $\mathfrak{g}$ of $X$:
\begin{equation}
\label{eq:6}
\begin{array}{lll} E_1 & = & a_{11}({x_3})\partial _{x_1}+a_{21}({x_3})\partial _{x_2}, \\
E_2& = & a_{12}({x_3})\partial _{x_1}+a_{22}({x_3})\partial _{x_2}, \\
 E_3 & = & \partial _{x_3}, \end{array}
\end{equation}
where we have denoted
\begin{equation}
\label{eq:exp(zA)}
 e^{{x_3}A}=\left(
\begin{array}{cr}
a_{11}({x_3}) & a_{12}({x_3}) \\
a_{21}({x_3}) & a_{22}({x_3})
\end{array}\right) .
\end{equation}


\begin{definition}
{\rm
Given a matrix $A\in {\mathcal M}_2(\R )$, we define the {\it canonical
left-invariant metric} on $\R^2\rtimes _A\R $ to be that one for which
the left-invariant basis $\{ E_1,E_2,E_3\} $ is orthonormal.
}
\end{definition}

For a general $X=\RA$, the isometry group of $X$ is generated by left translations in $X$ and the rotation $(x_1,x_2,x_3) \mapsto (-x_1,-x_2,x_3)$. However, for some matrices $A$ the isometry group is larger as occurs when $A=\left(
\begin{array}{cr}
1 & b \\
b & 1 \end{array}\right) $, $b\in [0,\infty)$, which gives a metric on $X=\RA$ of constant curvature $-1$ and any two of the underlying Lie groups are non-isomorphic; see~\cite{mpe11,mil2} for these results.

The Levi-Civita connection $\widehat{\nabla} $ for the
canonical left-invariant metric of $X=\R^2\rtimes _A\R $:
 {\large

\begin{equation} \label{eq:12}
\begin{array}{l|l|l}
\widehat{\nabla} _{E_1}E_1=a\, E_3 & \widehat{\nabla} _{E_1}E_2=\frac{b+c}{2}\, E_3 & \widehat{\nabla} _{E_1}E_3=-a\, E_1-\frac{b+c}{2}\, E_2 \\
\widehat{\nabla} _{E_2}E_1=\frac{b+c}{2}\, E_3 & \widehat{\nabla} _{E_2}E_2=d\, E_3 & \widehat{\nabla} _{E_2}E_3=-\frac{b+c}{2}\, E_1-d\, E_2 \\
\widehat{\nabla} _{E_3}E_1=\frac{c-b}{2}\, E_2 & \widehat{\nabla}
_{E_3}E_2=\frac{b-c}{2}\, E_1 & \widehat{\nabla} _{E_3}E_3=0.
\end{array}
\end{equation}
}

Finally we remark that the mean curvature of each leaf of the foliation ${\mathcal F}=
\{ \R^2\rtimes _A\{ {x_3}\} \mid {x_3}\in \R \} $ with respect to the unit
normal vector field $E_3$ is the constant $H=\mbox{trace}(A)/2$; this follows immediately from the above formulas for the connection and the fact that $E_1,E_2$ are tangent to the leaves of this foliation. Since we are assuming in this paper that $\mbox{trace}(A)\geq 0$, the planes in this foliation are mean convex in the direction of $E_3$.


\section{The proof of Theorem~\ref{main}.} \label{sec:halfspace}
Throughout this section, $X=\RA$, where $A=\left(
\begin{array}{cc}
a & b\\
c & d\\
 \end{array}\right).$  We will use coordinates $x_1,x_2,x_3$ for $\RA$ given in the previous section.   Recall that we assume that  $\mbox{trace}(A)=a+d\geq0.$  In particular, by the discussion at the end of the previous section, the mean curvature of the planes $\R^2 \rtimes_A \{z\}$ in $\RA$ have constant mean curvature $H_0(X)=\frac{a+d}{2}\geq0$ with respect to the unit normal field $\partial_{x_3}$.

We now come to a crucial lemma related to the geometry of a small neighborhood of the plane $\R^2 \rtimes_A \{0\}$ in $\RA$. Note that in these spaces the function $\phi$ described in the next lemma is just the inverse of the distance function to the normal subgroup $\R^2 \rtimes_A \{0\}$.

\begin{lemma} \label{lemma} Let $X=\RA$. There exists a constant $C_1>0$ such that for any $H$-surface $\Sigma$  in $X$ with $|H|\leq H_0(X)$ and satisfying $$0<x_3\leq C_1$$ on $\Sigma$, then the function $\phi \colon =\frac{1}{x_3}$ is subharmonic on $\Sigma$.
\end{lemma}

\begin{proof} We view $\Sigma$ as a conformal immersion $f=(x_1,x_2,x_3)\colon \Sigma \to X$ from a Riemann surface $\Sigma$ and let $z$ be a conformal coordinate. We define $A_1,A_2, A_3$ by $$f_z=A_1E_1+A_2E_2 +A_3E_3.$$
In particular, $A_3=x_{3z}$. The conformality of $f$ means that the equation \begin{equation} \label{eq:conformal} A_1^2+A_2^2+A_3^2=0
\end{equation} holds.

Since $f$ has mean curvature $H$, we have $${\widehat{\nabla}}_{f_{\overline{z}}}f_z= H(|A_1|^2+|A_2|^2+|A_3|^2)N,$$ where $N$ denotes the unit normal vector. From this and \eqref{eq:12}, we get $$A_{3\overline{z}}=\langle {\widehat{\nabla}}_{f_{\overline{z}}}E_3, f_z\rangle +\langle E_3,\widehat{\nabla}_{f_{\overline{z}}}f_z\rangle$$
$$=-a|A_1|^2-d|A_2|^2 -\frac{b+c}{2}(A_1\overline{A}_2+ \overline{A}_1 A_2)+HN_3(|A_1|^2+|A_2|^2+|A_3|^2),$$ where $N_3=\langle N, E_3\rangle.$

Since
$$\phi_{z\overline{z}}=-\left(\frac{A_3}{x_3^2}\right)_{\overline{z}}=\frac{2|A_3|^2-x_3A_{3\overline{z}}}{x_3^3},$$
we get  $$ x_3^3\phi_{z\overline{z}}=(2-HN_3x_3)|A_3|^2+x_3\left( \frac{a+d}{2}-HN_3\right)(|A_1|^2+|A_2|^2)$$ $$+x_3\frac{a-d}{2}(|A_1|^2-|A_2|^2)+ x_3\frac{b+c}{2}(A_1\overline{A}_2+\overline{A}_1 A_2).$$

The second term in the right-hand expression is nonnegative, since $x_3>0$, $|H|\leq H_0(X)=\frac{a+d}{2}$ and $-1\leq N_3 \leq 1$. On the other hand, by \eqref{eq:conformal} we have $$|A_3|^4=|A_1^2+A_2^2|=|A_1|^4+|A_2|^4+2  \mbox{Re}(A_1^2 \overline{A}_2^2)$$ $$\geq |A_1|^4+|A_2|^4-2|A_1|^2|A_2|^2 =(|A_1|^2-|A_2|^2)^2,$$ and

$$|A_3|^4\geq 2|A_1|^2|A_2|^2+2  \mbox{Re}(A_1^2\overline{A}_2^2)=(A_1\overline{A}_2+\overline{A}_1 A_2)^2,$$

From this we obtain that $$x^3_3\phi_{z\overline{z}}\geq \left (2-x_3\left( |H|+\frac{|a-d|}{2}+\frac{|b+c|}{2}\right)\right)|A_3|^2.$$
Consequently, we get $\phi_{z\overline{z}}\geq 0$, when $$0<x_3\leq \frac{2}{|H|+\frac{|a-d|}{2}+\frac{|b+c|}{2}}=\colon C_1.$$ This completes the proof of the lemma.
\end{proof}

We now apply the above lemma to prove Theorem~\ref{main}.  We first consider the case of item~\ref{it1} of the theorem.  Suppose
that $F$ has mean curvature 0.  In this case the arguments given in the proofs of Theorems~1.3 and 1.5 in~\cite{dmr1} with
Lemma~\ref{lemma} replacing the roles of the respective  Lemmas~2.1 and  3.1 in~\cite{dmr1}, prove that item~\ref{it1} of the theorem holds.
 So assume now that $H_0(X)>0$.

We first construct for a sufficiently small $\varepsilon>0$ and for each $R>1$ a certain piecewise smooth domain
$$\Omega(R)\subset W(\varepsilon)=\{(x_1,x_2,x_3)\mid 0\leq x_3\leq \varepsilon\}.$$ Let
$$S(\delta,R)=\{(x_1,x_2,x_3) \mid x_1^2+x_2^2=R^2,x_3=\delta\}.$$ Assume $\varepsilon$ is chosen sufficiently small so that there exists
a compact embedded annulus $\cA\subset W(\varepsilon)$ satisfying:
\begin{enumerate}
\item $\partial\cA=\cA\cap \partial W(\varepsilon)=S(0,1)\cup S(\varepsilon,1)$,
\item $\cA$ has constant mean curvature $H_0(X)$ and the mean curvature vector of $\cA$ points out of the bounded component of $W(\varepsilon)-\cA$.
\end{enumerate}
For each $R>1$, define $\Omega(R)$ to be the closure of the nonsimply-connected bounded domain of
$$W(\varepsilon)-\left(\cA\cup\bigcup_{g\in S(0,R)}g\cA\right).$$ The boundary of $\Omega(R)$ is piecewise smooth and consists of $\cB_1(R)=\Omega(R)\cap F$,
$\cA$, an analytic surface $\cC$ and   $\cB_2(R)=\Omega(R)\cap  F(\ve)$, where $F(\ve)$ is the left translate $(0,0,\varepsilon) F$ of $F$. With respect to the inward pointing normal to
$\Omega(R)$, the surfaces $\cB_1(R)$ and $\cA$ have mean curvature $H_0(X)$, $\cC$ has mean curvature greater than or equal to $H_0(X)$ and $\cB_2(R)$ has mean curvature $-H_0(X)$.

We write
$\partial\Omega(R)=\Sigma_1\cup\Sigma_2$ where $\Sigma_1\subset\partial\Omega(R)$ is the compact piecewise
smooth annulus with $\partial\Sigma_1=S(0,R)\cup S(\varepsilon,1)$ which  contains $\cA$ and  $\Sigma_2\subset\partial\Omega(R)$ is the closure of
$\partial\Omega(R)-\Sigma_1$ in $\partial\Omega(R)$; we orient $\Sigma_1$ and $\Sigma_2$ by the unit normal pointing into $\Omega(R)$.
We associate to the pair $(S(0,R),S(\varepsilon,1))$ the functional
$$T_1=\mathrm{Area}+2H_0(X)\cdot\mathrm{Volume}$$ defined on oriented surfaces $\Delta$ with $\partial\Delta=S(0,R)\cup S(\varepsilon,1)$
and which are homologous to $\Sigma_1$ and where $\mathrm{Volume}$ refers to the volume of the subdomain of $\Omega(R)$ with boundary
$-\Delta\cup\Sigma_1$, where $-\Delta$ is the 2-cycle $\Delta $ with the opposite orientation.

By the arguments in the proof of
Theorem~2 in~\cite{rose4} and using the fact that $\partial\Omega(R)$ is a good barrier for minimizing the functional $T_1$, there exists
a smooth surface $\Delta_1(R)\subset\Omega(R)$ with $\partial\Delta_1=(S(0,R),S(\varepsilon,1))$ which minimizes $T_1$, and furthermore, $\Delta_1(R)$ is stable. Also by Theorem~2 in~\cite{rose4}, $\Delta_1(R)$ has mean
curvature $H_0(X)$ with respect to the outward pointing normal to the subdomain of $\Omega(R)$ with boundary $\Delta_1(R)\cup\Sigma_1$.

Arguing by contradiction, we now prove that item~\ref{it2} of the theorem holds. Suppose $H\in[0,H_0(X)]$ and $\cS\subset X$ is an $H$-surface
properly immersed in $X$ that is contained on the mean convex side of $F$. After translations of both $\cS$ and $F$, assume that
$F=\R^2\rtimes_A\{0\}$ and that the distance from $\cS$ to $F$ is zero. Recall  that $H_0(X)=\frac12\mathrm{trace}(A)>0$ and so,
by the maximum principle, the $x_3$-coordinate of $\cS$ is strictly positive (unless $\cS=F$, in which case the conclusion of item~\ref{it2} holds).

Since $\cS$ is proper in $X$ and disjoint from $F$, then $\cS$ is at a distance greater than some small $\varepsilon>0$ from the
disk $\{(x_1,x_2,x_3)\mid x_1^2+x_2^2\leq1,x_3=0\}$. After possibly choosing a  smaller $\ve$, we may assume that $S$ is disjoint from the compact annulus $\cA$ defined previously and bounded by $ S(0,1) \cup S(1,\ve)$. Let $\Delta_1(R)$ be a surface with $\partial\Delta_1(R)=S(0,R)\cup S(\varepsilon,1)$
which minimizes the functional $T_1$ associated to $(S(0,R),S(\varepsilon,1))$. Note that $\Delta_1(R)$ must be disjoint from $\cS$, since
otherwise we can translate $\Delta_1(R)$ by $(0,0,-t)$ for $t\in[0,\varepsilon]$ and there will be a last point where $(0,0,-t)\Delta_1(R)$
intersects $\cS$ and this point will be an interior point of $\Delta_1(R)$, contradicting the maximum principle.

We claim that there exists a constant $c>0$ such that the areas of the surfaces $\Delta_1(R)$ are bounded by $cR^2$. Indeed, we have
$$\mathrm{Area}(\Delta_1(R))\leq T_1(\Delta_1(R))\leq T_1(\Sigma_1)\leq\mathrm{Area}(\Sigma_1)+\mathrm{Volume}(\Omega(R)).$$
Since $W(\varepsilon)$ is quasi-isometric to a Euclidean slab, the area of $\Sigma_1$ and the volume of $\Omega(R)$ are quadratic in $R$.
This proves the claim.

By curvature estimates for stable $H$-surfaces in~\cite{rst1} ($\Delta_1(R)$ is stable and two-sided) and by local area estimates
($\Delta_1(R)$ minimizes the functional $T_1$), a subsequence of the surfaces $\Delta_1(n)$, $n\in\N-\{0,1\}$, converges smoothly
on compact domains of $W(\varepsilon)$ to a properly embedded $H_0(X)$-surface $\Delta_1(\infty)$ in $W(\varepsilon)$ with boundary $S(\ve,1)$.
The quadratic bound on the areas of the surfaces $\Delta_1(R)$ also implies that $\Delta_1(\infty)$ has intrinsic area growth which is at most quadratic. By the results described in~\cite{gri1}, the quadratic area growth of $\Delta_1(\infty)$  means it is a parabolic Riemannian surface.

An application of Lemma~\ref{lemma} then gives that the function $x_3$ must have the constant value of $\varepsilon$ on $\Delta_1(\infty)$ (note
that in order to apply the lemma we need to choose $\varepsilon<C_1/2$ and to translate $\Delta_1(\infty)$ by $(0,0,C_1/2)$, so that the function
$\phi$ is subharmonic and bounded, and thus cannot admit an interior maximum or a maximum at infinity unless it is constant).
But this implies that $\cS$ must intersect $\Delta_1(\infty)$ transversely at some point, which contradicts
the fact that $\cS$ does not intersect any of the surfaces $\Delta_1(R)$. This contradiction proves that item~\ref{it2} of the theorem holds.

We now argue by contradiction to prove item~\ref{it3} of the theorem. Since the proof if similar to that of item~\ref{it2}, we only give a
sketch of it. Suppose $H\in[H_0(X),+\infty)$ and $\cS\subset X$ is an $H$-surface
properly embedded in $X$ such that $F$ is contained on the mean convex side of $\cS$. After translations of both $\cS$ and $F$, assume that
$F=\R^2\rtimes_A\{\varepsilon\}$ and that the distance from $\cS$ to $F$ is zero. By the maximum principle, the $x_3$-coordinate of $\cS$ is strictly less than $\varepsilon$ (unless $\cS=F$).

Since $\cS$ is proper in $X$ and disjoint from $F$, then $\cS$ is at a distance greater than some small $\varepsilon>0$ from the
disk $\{(x_1,x_2,x_3) \mid x_1^2+x_2^2\leq1,x_3=\varepsilon\}$.
We write
$\partial\Omega(R)=\Sigma_1\cup\Sigma_2$, where $\Sigma_1\subset\partial\Omega(R)$ is the compact piecewise
smooth annulus with $\partial\Sigma_1=S(0,1)\cup S(\varepsilon,R)$ that contains $\cA$ and  $\Sigma_2\subset\partial\Omega(R)$ is the closure of
$\partial\Omega(R)-\Sigma_1$ in $\partial\Omega(R)$; we orient $\Sigma_1$ and $\Sigma_2$ by the unit normal pointing into $\Omega(R)$.
We associate to the pair $(S(0,1),S(\varepsilon,R))$ the functional
$$T_2=\mathrm{Area}+2H_0(X)\cdot\mathrm{Volume}$$ defined on oriented surfaces $\Delta$ with $\partial\Delta=S(0,1)\cup S(\varepsilon,R)$
and which are homologous to $\Sigma_2$ and where $\mathrm{Volume}$ refers to the volume of the subdomain of $\Omega(R)$ with boundary
$-\Delta\cup\Sigma_2$.

As in the case of item~\ref{it2}, $\partial\Omega(R)$ is a good barrier for minimizing the functional $T_2$, so there exists
a smooth surface $\Delta_2(R)\subset\Omega(R)$ with $\partial\Delta_2=S(0,1)\cup S(\ve,R)$ which minimizes $T_2$. Also by the proof of Theorem~2 in~\cite{rose4}, $\Delta_2(R)$ has mean
curvature $H_0(X)$ with respect to the outward pointing normal to the subdomain of $\Omega(R)$ with boundary $-\Delta_2(R)\cup\Sigma_2$.
As in the case of item~\ref{it2}, $\Delta_2(R)$ must be disjoint
 from $\cS$ (here we use the fact that the
mean curvature vectors of both surfaces point towards $F$ at the intersection point).

A subsequence of the surfaces $\Delta_2(n)$, $n\in\N-\{0,1\}$, converges smoothly
on compact domains of $W(\varepsilon)$ to a properly embedded $H_0(X)$-surface $\Delta_2(\infty)$ in $W(\varepsilon)$ with boundary $S(0,1)$, and
$\Delta_2(\infty)$ has at most quadratic area growth and thus has parabolic conformal type. An application of Lemma~\ref{lemma} then
shows that $x_3$ is constant on $\Delta_2(\infty)$ and the conclusion is the same as in the case of item~\ref{it2}. This final contradiction completes the proof of Theorem~\ref{main}.

\section{Partial results on the embedded Calabi-Yau problem.} \label{sec:CY}
In this section we will prove Theorem~\ref{CY}  by applying Theorem~\ref{main} and Theorem~\ref{CY2} by applying the results in \cite{dmr1}. The main techniques needed to prove this theorem arise from the study of the geometry and structure of embedded $H$-disks in a homogeneous three-manifold. This study when $H=0$ was initiated by Colding and Minicozzi in a series of papers \cite{cm25,cm23,cm35} and then generalized to the case $H>0$ by Meeks and Tinaglia \cite{mt8,mt7,mt1,mt9}.  The proofs of Theorems~\ref{CY} and~\ref{CY2} depend on several key related results which we now recall.

In~\cite{cm35}, Colding and Minicozzi proved a general one-sided curvature estimates for  embedded minimal disks, which was subsequently generalized by Meeks and Tinaglia~\cite{mt8,mt9} for $H>0$. An important application of this one-sided curvature estimate is to prove that whenever $\cS$ a complete embedded $H$-surface with positive injectivity radius in a complete Riemannian three-manifold $Y$, then on compact subsets $\Delta$ of $Y$, the second fundamental form of $\cS\cap \Delta$ is bounded; see~\cite{mr13} for this application of one-sided curvature estimates and other results in~\cite{cm35} when $H=0$ and~\cite{mt8} for the case $H>0$. It follows that the closure $\overline{\cS}$ of such an $H$-surface $\cS$ has a lamination type structure called a weak $H$-lamination. Basically  a weak $H$-lamination structure for $\overline{\cS}$ means  that $\overline{\cS}$ is the union of complete embedded $H$-surfaces in $Y$ with second fundamental form  bounded on compact domains of $Y$ and we refer the reader to~\cite{mpr19,mpr18} for further discussion on the notion of a weak $H$-lamination; the word {\em weak} is used because two leaves of this type of ``lamination" may intersect non-transversely along zero or one-dimensional subsets, with two such intersecting leaves at a point $p$ lying   locally on one side of each other at $p$.

Two-sided covers of limit leaves of a weak $H$-lamination $\overline{\cS}$ of $Y$ are stable~\cite{mpr18} and so, by the universal curvature estimates for two-sided stable $H$-surfaces in~\cite{rst1}, the leaves of the sub weak $H$-lamination ${\rm Lim}(\overline{\cS})$ of limit leaves of $\overline{\cS}$ have uniformly bounded second fundamental forms, where the bound only depends on a bound of the absolute sectional curvature of $Y$ in some neighborhood of ${\rm Lim}(\overline{\cS})$. An application of the one-sided curvature estimates for embedded $H$-disks implies that for $Y$ having bounded absolute curvature, there exists a small $\ve>0$, such that the norm of the second fundamental forms of the leaves of $\overline{\cS}$ are uniformly bounded in the $\ve$-neighborhood of ${\rm Lim}(\overline{\cS})$. Shortly we will apply this  curvature estimate result for $\overline{\cS}$ near ${\rm Lim}(\overline{\cS})$  in our proofs of Theorems~\ref{CY} and~\ref{CY2} but first we recall two key related results.

\begin{theorem}[Meeks-Perez, \,\cite{mpe11}]
\label{mp}
A complete embedded minimal surface $\cS$ of finite topology in a homogeneous three-manifold has positive injectivity radius. In particular, if $\cS$ is not proper, then the closure $\overline{\cS}$ has the structure of a minimal lamination with nonempty sublamination ${\rm Lim}(\overline{\cS})$ and in some small $\ve$-neighborhood $N(\ve)$ of ${\rm Lim}(\overline{\cS})$, $\cS \cap N(\ve)$ has bounded second fundamental form.
\end{theorem}

\begin{theorem}[Meeks-Tinaglia, \,\cite{mt1}] \label{mt} The closure of a complete embedded $H$-surface of finite topology in a homogeneous three-manifold $X$ has the structure of a weak $H$-lamination; furthermore, $H$-surfaces of finite topology in $X$ have bounded second fundamental form on compact subsets of $X$.

More generally, if $\cS$ is a complete embedded $H$-surface in a homogeneous three-manifold $X$, then the following properties are equivalent:
\begin{enumerate}
\item the closure $\overline{\cS}$ of $\cS$ has the structure of a weak $H$-lamination,
\item the injectivity radius function of  $\cS$ is bounded away from zero on compact subsets of $X$,
\item $\cS$  has bounded second fundamental form on compact subsets of $X$.
\end{enumerate}
\end{theorem}


We are now in a position to prove Theorem~\ref{CY}.

\begin{proof}[Proof of Theorem~\ref{CY}]

Fix a value $H \in[0,H_0(X)]$.
Let $F$ be a translate of a normal $\R^2$ subgroup of a simply-connected metric Lie group $X=\RA$. Note that $X$ is equipped with global $\rth=\R^2\times \R=\{(x_1,x_2,x_3) \mid x_1,x_2,x_3\in \R\}$ coordinates, which we will use in the following discussion.

Arguing by contradiction, suppose $\cS$ is a complete, connected embedded $H$-surface which is contained on a (the) mean convex side of $F$ but which is {\em not} a left translate of $F$. After left translating $\cS$, we may assume that the distance from $\cS$ to $F$ is zero. We let $W$ denote the closure of the mean convex component of $X-F$ containing $\cS$.

First suppose $\cS$ has positive injectivity radius. By Theorem~\ref{main}, $\cS$ is not properly embedded in $X$ and so Theorem~\ref{mt} implies $\overline{\cS}$ has the structure of a weak $H$-lamination of $X$ with ${\rm Lim}(\overline{\cS})\neq \O$.

Using the facts that $\overline{\cS}$ is a closed subset of $X$ and the weak $H$-lamination $\overline{\cS}$ has distance zero to $\partial W=F$,  the arguments in the proof of Theorem~\ref{main} using the maximum principle imply that the plane $F$ must be a limit leaf of $\overline{\cS}$. By our previous discussion, $\overline{\cS}$ has uniform curvature estimates in a fixed size closed $\ve$-neighborhood $N(\ve)$ of $F$ in $W$. It follows that for $\ve$ sufficiently small, each component $\Delta$ of $\cS\cap N(\ve)$ is a multigraph over its $x_3$-projection to $F$ and without loss of generality we will assume that $F=\R^2\rtimes_A\{0\} \subset \RA =X$.

Elementary covering space arguments imply that $\Delta$ is actually a graph over its $x_3$-projection to $F$, see Lemma~1.1 in~\cite{mr8} for this argument. Since such a graph with boundary is proper in the closed set $N(\ve)$ with boundary in the plane $(0,0,\ve)F\subset \partial N(\ve)$, we can apply the arguments in the proof of Theorem~\ref{main} to obtain a contradiction. This contradiction proves the theorem when $\cS$ has positive injectivity radius.
By Theorem~\ref{mp}, if $\cS$ has finite topology and $H=0$, then $\cS$ has positive injectivity radius. Moreover, by
Theorem~\ref{mt}, if $\cS$ has finite topology then its injectivity radius function is bounded away from zero on compact subsets of $X$.
Hence, to complete the proof of the theorem, it suffices to assume now that $H>0$ and to prove the theorem in the case of item~2.
So, suppose now that $H \in(0,H_0(X)]$ and $\cS$  has  positive injectivity radius function on compact domains of $X$.

Since Theorem~\ref{mt} implies that any complete, embedded finite topology $H$-surface in $X$ has injectivity radius function bounded away from zero on compact subsets of $X$, we only have  the case of item~2 of the  theorem left to prove. So, suppose that $H \in(0,H_0(X)]$ and $\cS$ is a complete, connected embedded $H$-surface  in a (the) mean convex component of $X-F$, $\cS$ has  positive injectivity radius function on compact domains of $X$ and $\cS$ is not a left translate of $F$. Let $W$ be the closure of this mean convex component.

Before continuing the proof,  we demonstrate a general properness property for a complete, connected embedded $H$-surface (including the case $H=0$) $\Sigma$ satisfying:
\bit \item $\Sigma$ is contained in a mean convex component of $X-F$; \item  $\Sigma$ has injectivity radius function bounded away from zero on compact subsets of $X$;  \item $\Sigma$  is not a translate of $F$; \item   $\Sigma $ is not contained in a smaller halfspace of $X-F$. \eit By Theorem~\ref{mt}, the closure $\overline{\Sigma}$ has the structure of a weak $H$-lamination  with ${\rm Lim}(\overline{\Sigma})\neq \O$.  Since any leaf $L$ of ${\rm Lim}(\overline{\Sigma})$ is stable, then such an $L$ has bounded second fundamental form and our previous arguments imply $L$ is a translate of $F$.  In particular, by the connectedness of $\Sigma$ either $\Sigma$ is proper in $X-F$ or there exists a left translate $F' \subset (X-F)$ of $F$ and $\Sigma$ is proper in the open slab with boundary $F\cup F'$.  Since each component of $X-F$ (resp. $X-[F\cup F']$) is simply-connected, then elementary separation properties imply that \begin{quote} {\bf $\Sigma$ separates the halfspace or slab of $X$ in which it is contained as a proper surface. }\end{quote}

We now return to the proof of the case when $H>0$ and $\cS$ is a complete, connected embedded $H$-surface  in the mean convex component of $X-F$, $\cS$ has  injectivity radius function bounded away from zero on compact domains of $X$ and $\cS$ is not a left translate of $F$; note that Theorem~\ref{mt} implies that this case for $\cS$ includes the case where $\cS$ has finite topology. By the discussion in the previous paragraph,  $F$ is a limit leaf of the weak $H$-lamination $\overline{\cS}$ and $\cS$ separates a halfspace or slab of $X$ in which it is contained as a proper surface.  It follows that close to any compact disk $D$ of $F$, the mean curvature vectors of successive disks contained in leaves of the weak $H$-lamination  $\overline{\cS}$, and which are converging to $D$, have non-zero mean curvature vectors with alternating positive and negative $x_3$-coordinates.   But the mean curvatures of this sequence of disks converge to the {\bf non-zero} mean curvature $D$, which gives a contradiction.  Theorem~\ref{CY} now follows from this contradiction. \end{proof}

We now outline the proof of Theorem~\ref{CY2} stated in the Introduction as the arguments in its proof are similar to the proofs of the two statements in Theorem~\ref{CY}. Let $\cS$ be a complete embedded minimal surface in $\nil$ which lies on one side of an entire minimal graph $F$. Since $F$ is stable, it satisfies a uniform curvature estimate. On the other hand, we may assume that  $\cS$ intersects the vertical $\ve$-translate of $F$ for every
$\ve$ small enough and so, by one-sided curvature estimates for embedded minimal disks, $\cS$ has uniformly bounded curvature in a fixed size small regular neighborhood of $F$. It follows from our discussion of the similar cases in the proof of Theorem~\ref{CY} that for some small $\ve>0$, the intersection of the closed $\ve$-neighborhood ${N(\ve)}$ of $F$ with $\cS$ consists of components which are proper normal graphs over proper subdomains of $F$. Then, by our previous arguments in the proof of Theorem~\ref{CY} and of Theorem~1.4 in~\cite{dmr1} in the case that $\cS$ was proper in $\nil$, we obtain a contradiction; this is the same argument that we use to prove Theorem~\ref{CY} when $\cS$ has positive injectivity radius. This concludes our sketch of the proof of Theorem~\ref{CY2}.

\begin{remark}{\rm
Let $M$ be a complete, simply-connected parabolic Riemann surface with bounded Gaussian curvature. Then the techniques developed in~\cite{dmr1} can be used to prove that a properly immersed minimal surface in $M\times\R$ that does not intersect $M\times\{0\}$ is some slice $M\times\{t\}$; a special case was obtained in~\cite{rose2}. For this we use the fact the third coordinate is harmonic on minimal surfaces, that vertical translations are isometries and that $M\times\R$ has bounded geometry to have uniform curvature estimates.
Then, the techniques developed in this section show that a complete embedded minimal surface of finite topology or of positive injectivity radius in $M\times\R$ that does not intersect $M\times\{0\}$ is some slice $M\times\{t\}$.

Let $M$ be as in the previous paragraph. Using the homological invariance of flux of the harmonic $t$-coordinate function of a complete embedded minimal surface $\cS\subset M\times \R$ along simple closed curves in $\cS$, it is can be proved that if $M$ has finite topology, then $M$ has positive injectivity radius; see \cite{mpe11,mr13} for the necessary techniques needed to prove this positive injectivity radius property. Hence, by the discussion in the previous paragraph, a complete embedded minimal surface $\cS\subset [(M\times \R) -(M\times \{0\})]$ of finite topology is equal to  $M\times\{t\}$ for some $t$. Related to this discussion and to Conjecture~\ref{con:CY} in the Introduction, we conjecture: If $M$ is a complete, simply-connected parabolic Riemannian surface and $\cS\subset M\times \R$ is complete embedded minimal surface of finite topology, then $\cS$  is  proper in $M\times \R$. }
\end{remark}

\bibliographystyle{plain}
\bibliography{dmr}

\end{document}